\documentclass[12pt]{article}
\usepackage{latexsym}
\usepackage[all]{xy}
\usepackage{amsmath}
\usepackage{amssymb}
\usepackage{amsfonts}

\newtheorem{thm}{Th\'eor\`em}[section]
\newtheorem{prop}[thm]{Proposition}
\newtheorem{lem}[thm]{Lemme}
\newtheorem{sublem}[thm]{Sous-lemme}
\newtheorem{df}[thm]{D\'efinition}
\newtheorem{cor}[thm]{Corollaire}

\begin{document}

\title{\textbf{Vers une axiomatisation de 
la th\'eorie des 
cat\'egories sup\'erieures}}
\bigskip
\bigskip

\author{\bigskip\\
Bertrand To\"en \\
\small{Laboratoire Emile Picard}\\
\small{UMR CNRS 5580} \\
\small{Universit\'{e} Paul Sabatier, Toulouse}\\
\small{France}}

\date{March 2005}

\maketitle

\begin{abstract}
On d\'efinit une notion axiomatique de \emph{th\'eorie
de $(1,\infty)$-cat\'egories}. On d\'emontre qu'une telle th\'eorie  
est  unique \`a \'equivalence (essentiellement unique) pr\`es. 
\end{abstract}

\medskip

\tableofcontents

\bigskip

\section{Introduction}

La th\'eorie des cat\'egories sup\'erieures a 
connu ces derniers temps un regain d'int\'er\^et, et 
il existe aujourd'hui au moins une dizaine
de d\'efinitions diff\'erentes de $n$-cat\'egorie (voir 
par exemple \cite{le}). Une question fondamentale
est bien entendu de les comparer. De prime abord cette question 
ne parra\^it pas tout \`a fait \'evidente \`a poser
correctement, car il est bien connu que les
$n$-cat\'egories doivent elles-m\^eme former
une $(n+1)$-cat\'egorie, et de ce fait comparer
deux th\'eories de $n$-cat\'egories demanderait
\`a savoir d\'ej\`a comparer les notions
de $(n+1)$-cat\'egories
de ces m\^emes th\'eories. Cette \emph{poursuite des cat\'egories}
ne semble mener nulle part. 

Cependant, il semble maintenant admis qu'avant 
m\^eme de former une $(n+1)$-cat\'egorie, les $n$-cat\'egories
doivent former une cat\'egorie de mod\`eles, ou tout
au moins une certaine sous-cat\'egorie raisonable. 
Ce point de vue, en plus d'\^etre particuli\`erement 
efficace pour d\'evelopper des th\'eories de $n$-cat\'egories
(voir par exemple \cite{jo,hs,pe}), 
poss\`ede aussi l'avantage de permettre 
de poser clairement la question de la
comparaison entre deux th\'eories: on dira que deux
th\'eories de $n$-cat\'egories sont \'equivalentes si
les cat\'egories de mod\`eles correspondantes
sont \'equivalentes au sens de Quillen. 
Suivant cette id\'ee, P. May a r\'ecemment propos\'e
un vaste programme de recherche
consistant \`a construire 
des \'equivalences de Quillen entre les
diff\'erentes th\'eories
de $n$-cat\'egories. 

La but de cet article et d'amorcer ce travail de
comparaison en essayant de d\'egager des
axiomes caract\'erisant, \`a \'equivalence
de Quillen pr\`es, "la" cat\'egorie
de mod\`eles des $n$-cat\'egories \footnote{L'id\'ee que ces
axiomes puissent exister n'est pas nouvelle (voir
\cite{si}).}. \\

\begin{center} \textit{$(1,\infty)$-Cat\'egories} \end{center}

Dans cet travail, nous ne nous int\'eresserons pas \`a 
la th\'eorie des $n$-cat\'egories pour 
$n$ quelconque (pouvant \^etre infini), mais au cas
particulier des \emph{$(1,\infty)$-cat\'egories}. 
De fa\c{c}on intuitive, les $(1,\infty)$-cat\'egories sont 
des $\infty$-cat\'egories pour les quelles 
les $\infty$-cat\'egories de morphismes entre deux objets
fix\'es sont des $\infty$-groupoides. Autrement dit,
une $\infty$-cat\'egorie est une $(1,\infty)$-cat\'egorie
si tous ses $i$-morphismes sont inversibles
pour $i>1$. Il se trouve que de tr\`es nombreuses
$\infty$-cat\'egories que l'on rencontre dans la nature
sont des $(1,\infty)$-cat\'egories, et ceci car toute $\infty$-cat\'egorie
obtenue \`a partir d'une (1-)cat\'egorie en inversant formellement
un certain ensemble de fl\`eches est de ce type.  
La th\'eorie des $(1,\infty)$-cat\'egories et de ce fait riche
d'exemples. De plus, il existe d\'ej\`a plusieurs
notions mod\`elisant les $(1,\infty)$-cat\'egories et qui forment
des cat\'egories de mod\`eles, comme par exemple les
cat\'egories simplicialement enrichies, les $1$-cat\'egories
de Segal, les quasi-cat\'egories ou encore les
espaces de Segal complets (voir
\cite{be,hs,jo,re,pe}). Ces diff\'erentes th\'eories ne sont 
pour le moment pas compar\'ees, bien que conjecturalement
\'equivalentes. \\

\begin{center} \textit{Th\'eories de $(1,\infty)$-cat\'egories} \end{center}
 
La notion fondamentale de ce travail est celle
de \emph{th\'eorie de $(1,\infty)$-cat\'egories}. Une telle th\'eorie
est par d\'efinition la donn\'ee d'une cat\'egorie de mod\`eles
$M$ et d'un objet co-simplicial $C$ dans $M$, tel que certains
axiomes $A1-A7$ soient v\'erifi\'es (voir
Def. \ref{dgiraud} pour plus de d\'etails). 
Nous ne reproduirons pas ces axiomes dans cette introduction, mais
signalons qu'ils sont des analogues des axiomes
de type Giraud caract\'erisant le topos des ensembles\footnote{Le topos
des ensembles est caract\'eris\'e par les axiomes de Giraud
plus l'axiome affirmant que l'objet final $*$ est non vide et g\'en\'erateur.}. 
La
nouveaut\'e dans le contexte des
$(1,\infty)$-cat\'egories est l'existence 
de l'objet co-simplicial $C$, qui joue le r\^ole
de la cat\'egorie $I$ ayant deux objets et un unique morphisme
entre eux munie de sa structure de co-cat\'egorie. Cet
objet permet de d\'efinir la notion de quotient d'une 
\emph{relation cat\'egorique} dans $M$, qui est une 
extension de celle de quotient par des objets en groupoides.
Une fois cette notion d\'egag\'ee, les axiomes $A1-A5$
sont essentiellement les axiomes de Giraud
(au sens des topoi de mod\`eles de \cite{hagI}) mais
o\`u l'effectivit\'e des relations d'\'equivalences
est remplac\'ee par une effectivit\'e des
relations cat\'egoriques. Les deux derniers
axiomes A6 et A7 sont sp\'ecifiques \`a notre cas, et sont
des conditions de g\'en\'eration et non d\'eg\'en\'erescence.

Un exemple de th\'eorie des $(1,\infty)$-cat\'egories est 
fournit par la cat\'egorie de mod\`eles $CSS$ 
des espaces de Segal complets de C. Rezk (voir \cite{re}).
Le th\'eor\`eme principal de ce travail est le suivant.

\begin{thm}\label{ti}
Toute th\'eorie de $(1,\infty)$-cat\'egories
est \'equivalente (au sens de Quillen) \`a $CSS$.
\end{thm}

\bigskip

Dans le dernier paragraphe nous pr\'esentons aussi 
un argument qui montre que la th\'eorie
de $(1,\infty)$-cat\'egories $CSS$ ne poss\`ede pas
d'auto-\'equivalences. Pr\'ecisemment, nous montrons que
la localis\'ee simpliciale $LCSS$, en tant qu'objet de
la cat\'egorie de mod\`eles des $S$-cat\'egories ne poss\`edent aucune
auto-\'equivalences fixant la
sous-cat\'egorie $\Delta \hookrightarrow LCSS$. On d\'eduit de ceci
que le seules auto-\'equivalences de $LCSS$ sont celles
induites par des auto-\'equivalences de la cat\'egorie $\Delta$, qui 
sont r\'eduites au groupe $\mathbb{Z}/2$ (voir Thm. \ref{t2}). 
Ceci compl\`ete
le th\'eor\`eme \ref{ti} et implique que
si $M$ est une th\'eorie de $(1,\infty)$-cat\'egories
alors il existe essentiellement une unique \'equivalence
entre les th\'eories homotopiques associ\'ees \`a $M$ et \`a $CSS$. \\

\bigskip

\begin{center} \textit{Relations avec d'autres travaux}
\end{center}

Une premi\`ere source d'inspiration pour ce travail
a \'et\'ee l'article \cite{si}, dans lequel
des propri\'et\'es conjecturalement caract\'eristiques
de la th\'eorie des $n$-cat\'egories sont \'enonc\'ees.
Cependant, nos axiomes sont tr\`es diff\'erents
que ceux de \cite{si}, et on ne peut raisonablement 
consid\'erer notre th\'eor\`eme
\ref{ti} comme une r\'eponse possible
\`a la conjecture principale (Conjecture 3) de \cite{si}. 

Ma seconde source d'inspiration f\^ut 
les quelques discussions que j'ai pu avoir avec G. Vezzosi autour
des questions \emph{Qu'est-ce qu'une topologie sur une
$n$-cat\'egorie ?} et \emph{Qu'est-ce qu'un $n$-topos ?} 
Ce sont ces discussions qui m'ont donn\'e une id\'ee de
la forme que les axiomes de Giraud devraient avoir
dans le contexte des cat\'egories sup\'erieures. Ainsi, les
axiomes $A1-A5$ de la d\'efinition \ref{dgiraud}
donnent une petite id\'ee sur ce que pourrait \^etre
un \emph{$(2,\infty)$-topos} (en supposant par
convention que les \emph{$(1,\infty)$-topoi} sont
les topoi de mod\`eles de \cite{hagI}), bien 
que la forme actuelle de ces axiomes doivent 
probablement \^etre l\'eg\`erement
modifi\'ee afin d'\^etre plus facilement \emph{internalisables}. 

Enfin, pour ce qui est du probl\`eme de la comparaison
des diff\'erentes th\'eories de $n$-cat\'egories,
je n'ai pas v\'erifi\'e en d\'etails que les
cat\'egories de mod\`eles 
des cat\'egories simpliciales, des quasi-cat\'egories
et des $1$-cat\'egories de Segal \'etaient
des th\'eories de $(1,\infty)$-cat\'egories. Cela ne semble
pas tr\`es difficile mais n'est pas tout \`a fait \'evident non plus, 
et demanderait un article en soit
(le point d\'elicat semble \^etre l'axiome $A5$). Quand \`a l'extension 
de la notion de th\'eorie de $(1,\infty)$-cat\'egories au cas
des $(n,\infty)$-cat\'egories (i.e. des $\infty$-cat\'egories
dont les $i$-morphismes sont inversibles d\`es que $i>n$), 
il semble clair que cela doit passer tout d'abord par une
notion de $n$-intervalle, qui devrait \^etre un objet
multi-co-simplicial $\Delta^{n} \longrightarrow M$
satisfaisant certaines propri\'et\'es (comme
par exemple \^etre une $n$-co-cat\'egorie faible). Une fois
cette notion d\'egag\'ee la g\'en\'eralisation 
des axiomes $A1-A7$ \`a ce contexte me parrait un exercice
acad\'emique. Il est donc tout \`a fait raisonable de penser
g\'en\'eraliser le th\'eor\`eme 
\ref{ti} au contexte des $(n,\infty)$-cat\'egories
sans trop de probl\`emes. Cependant, montrer que 
les notions de $n$-cat\'egories d\'ecrites dans
\cite{le} v\'erifient les axiomes des th\'eories de $(n,\infty)$-cat\'egories 
est une autre histoire. \\

\bigskip

\begin{center} \textit{Remerciements}
\end{center}

Je voudrais remercier J. Lurie et 
G. Vezzosi pour plusieurs discussions
autour de la notion de topos sup\'erieurs, qui 
ont inspir\'e en partie ce travail.

Je tiens aussi a remercier le rapporteur 
pour ces remarques, et pour m'avoir signaler
un argument incomplet dans la preuve
du th\'eor\`eme principal. \\

\bigskip

\bigskip

\textbf{Conventions et notations:}
Nous utiliserons l'expression \emph{cmf} pour signifier
\emph{cat\'egorie de mod\`eles ferm\'ee}, et ce au
sens de \cite{ho}. La cmf des ensembles simpliciaux
sera not\'ee $SEns$.
Pour une cmf $M$, nous noterons
$Ho(M)$ sa cat\'egorie homotopique, et l'ensemble des morphismes
entre deux objets $x$ et $y$ dans $Ho(M)$ sera
not\'e $[x,y]_{M}$ (ou encore $[x,y]$ s'il n'y a
pas d'ambiguit\'e). Les \emph{espaces de morphismes},
tel que definis dans \cite[\S 5]{ho},
entre deux objets $x$ et $y$ dans $M$ seront
not\'es $Map_{M}(x,y)$ (ou encore $Map(x,y)$ s'il n'y a
pas d'ambiguit\'e), et toujours consid\'er\'es comme
objets de la cat\'egorie homotopique $Ho(SEns)$.
Les objets $Map(x,y)$ font de $Ho(M)$ une cat\'egorie
enrichie dans $Ho(SEns)$. Cet enrichissement est de plus
ferm\'e, au sens o\`u il existe pour $K\in Ho(SEns)$ et
$x,y\in Ho(M)$ des objets $K\otimes^{\mathbb{L}}x\in Ho(M)$ et
$y^{\mathbb{R}K}$ et des isomorphismes fonctoriels
$$Map_{SEns}(K,Map_{M}(x,y))\simeq
Map_{M}(K\otimes^{\mathbb{L}}x,y)\simeq Map_{M}(x,y^{\mathbb{R}K}).$$
Les produits fibr\'es homotopiques dans $M$ seront
not\'es $x\times^{h}_{z}y$ (en particulier
nous noterons $x\times^{h}y$ le produit homotopique).
De m\^eme, les sommes amalgam\'ees
homotopiques seront not\'ees $x\coprod^{\mathbb{L}}_{z}y$.

Enfin, nous dirons que deux cmf sont \emph{Quillen-\'equivalentes}
s'il existe une chaine d'\'equivalences de Quillen entre
les deux (dans des directions arbitraires).

\section{Cat\'egories de mod\`eles internes}

Commen\c{c}ons par rappeler l'existence de la notion de
cmf monoidale. Pour cela nous renvoyons \`a \cite[\S 4]{ho}.

\begin{df}\label{d1}
Une cmf est \emph{interne} si le produit direct en
fait une cmf monoidale.
\end{df}

En clair, une cmf $M$ est interne si elle v\'erifie
les deux conditions suivantes.
\begin{itemize}
\item La cat\'egorie $M$ est cart\'esiennement
close. En d'autres termes, pour deux
objets $x$ et $y$ dans $M$, le foncteur
$z \mapsto Hom(z\times x,y)$ est repr\'esentable
par un objet $\underline{Hom}_{M}(x,y)\in M$.
\item Pour toute paire de cofibrations
$$u : x \longrightarrow y \qquad  v : a \longrightarrow b,$$
le morphisme
$$u\Box v : x\times b \coprod_{x\times a}y\times a \longrightarrow
y\times b$$
est une cofibration, qui est de plus une cofibration triviale
si $u$ ou $v$ est une cofibration triviale.
\end{itemize}

Le lecteur remarquera que nous n'avons pas mention\'e 
l'axiome habituel des cmf monoidales concernant
l'unit\'e (voir \cite[\S 4]{ho}). En effet, il se d\'eduit
du deuxi\`eme axiome ci-dessus de la fa\c{c}on suivante.
Soit $x$ un objet cofibrant, et $Q(*) \rightarrow *$ un remplacement
cofibrant de $*$. Soit $x\rightarrow x'$ un remplacement
fibrant de $x$. Comme $Q(*)$ est cofibrant, le morphisme induit
$$Q(*)\times x \longrightarrow Q(*)\times x'$$
est une cofibration triviale et en particulier une \'equivalence. 
De plus, $x'$ \'etant fibrant, le morphisme induit
$$Q(*) \times x' \longrightarrow *\times x'\simeq x'$$
est une fibration triviale, et donc une \'equivalence. 
Le diagramme commutatif suivant
$$\xymatrix{
Q(*)\times  x \ar[r] \ar[d] & x \ar[d] \\
Q(*) \times x'\ar[r] & x'}$$
montre alors que $Q(*)\times x\longrightarrow x$ est une
\'equivalence. \\

Lorsque $M$ est une cmf interne, les
$Hom$ internes d\'eriv\'es sont des
$Hom$ internes pour la cat\'egorie homotopique $Ho(M)$
(voir \cite[\S 4]{ho}). De plus, ces $Hom$ internes
sont compatibles avec l'enrichissement naturel
de $Ho(M)$ dans $Ho(SEns)$, au sens qu'il exite
des isomorphismes fonctoriels dans $Ho(SEns)$
$$Map(z\times^{h}x,y)\simeq Map(z,\mathbb{R}\underline{Hom}_{M}(x,y)).$$
Une cons\'equence (presque) imm\'ediate de ceci est que pour tout
diagramme d'objets de $M$, $x_{\bullet} : I \longrightarrow M$,
le morphisme naturel
$$Hocolim_{i\in I}(x_{i}\times^{h}y) \longrightarrow
(Hocolim_{i\in I}x_{i})\times^{h}y$$
est un isomorphisme dans $Ho(M)$.

\begin{df}\label{d2}
Une cmf $M$ est \emph{faiblement interne}
si pour tout $y\in M$, et tout diagramme
$x_{\bullet} : I \longrightarrow M$, le morphisme naturel
$$Hocolim_{i\in I}(x_{i}\times^{h}y) \longrightarrow
(Hocolim_{i\in I}x_{i})\times^{h}y$$
est un isomorphisme dans $Ho(M)$.
\end{df}

Bien entendu, une cmf interne est faiblement interne.
La r\'eciproque n'est pas vrai, mais nous avons
tout de m\^eme la proposition
suivante.

\begin{prop}\label{p1}
\begin{enumerate}
\item Toute cmf Quillen-\'equivalente \`a une cmf
interne est faiblement interne.
\item Soit $M$ une cmf combinatoire au sens
de \cite{du}. Alors, $M$ est faiblement interne si et seulement si
elle est Quillen-\'equivalente \`a une cmf simpliciale, interne, et
dans laquelle tous les objets sont cofibrants.
\end{enumerate}
\end{prop}

\textit{Preuve:} $(1)$ Ceci est imm\'ediat car les \'equivalences
de Quillen pr\'eservent les limites et colimites homotopiques. \\

$(2)$ Soit $M$ une cmf combinatoire et faiblement interne. Il nous
faut montrer qu'elle est Quillen-\'equivalente \`a une cmf interne,
simpliciale, et 
dans laquelle tous les objets sont cofibrants. Tout d'abord, en
appliquant \cite{du} on peut supposer que $M$ est une cmf
simpliciale dans laquelle tous les objets sont cofibrants. Nous
noterons alors $\underline{Hom}$ les $Hom$ simpliciaux de $M$.

Soit $C\subset M$ une sous cat\'egorie pleine
qui est petite et dont l'ensemble des objets est un 
ensemble de repr\'esentant
pour les classes d'isomorphie d'objets
$\lambda$-petits pour un cardinal r\'egulier
$\lambda$ comme dans \cite[Def. 2.2]{du}. Quitte \`a choisir
$\lambda$ assez grand, on peut supposer que
$C$ est stable par produits directs ainsi que par le
foncteur de remplacement fibrant $R$. On note alors
$C^{f}$ la sous-cat\'egorie pleine de $C$ form\'ee
des objets fibrants. On consid\`ere alors le plongement
de Yoneda simplicial restreint
$$\begin{array}{cccc}
\underline{h} : & M & \longrightarrow & SPr(C^{f}) \\
 & x & \mapsto & \left( \underline{h}_{x}:=\underline{Hom}(-,x)\right) ,
\end{array}$$
o\`u $SPr(C^{f})$ est la cmf des pr\'efaisceaux simpliciaux
sur $C^{f}$ (on utilise ici la structure projective pour
laquelle les fibrations et \'equivalences sont d\'efinis
niveaux par niveaux).
Le foncteur $\underline{h}$ poss\`ede un adjoint \`a gauche
$$Re : SPr(C^{f}) \longrightarrow M$$
qui en fait un foncteur de Quillen \`a droite. Tout comme dans
\cite[Thm. 4.2.3]{hagI}, on voit que le foncteur d\'eriv\'e
$$\mathbb{R}\underline{h} : Ho(M) \longrightarrow
Ho(SPr(C^{f}))$$
est pleinement fid\`ele. On sait donc par \cite[Prop. 3.2]{du} qu'il
existe un petit  ensemble  de morphismes $S$ dans
$SPr(C^{f})$ tel que l'adjonction de Quillen
$(Re,\underline{h})$ induise une \'equivalence de Quillen
entre $M$ est la localis\'ee de Bousfield \`a gauche
$L_{S}SPr(C^{f})$. On peut donc supposer que
$M=L_{S}SPr(C^{f})$.

La cat\'egorie $SPr(C^{f})$ poss\`ede aussi une autre structure
de cmf pour laquelle les cofibrations et les \'equivalences
sont d\'efinis niveaux par niveaux, et en particulier
pour laquelle tout objet est cofibrant. On note cette structure
$SPr_{inj}(C^{f})$, que l'on voit tout de suite \^etre
une cmf interne (car $SEns$ est interne).
Le foncteur identit\'e induit une \'equivalence
de Quillen entre $L_{S}SPr(C^{f})$ et
$L_{S}SPr_{inj}(C^{f})$ et nous sommes donc ramen\'e \`a d\'emontrer
que $L_{S}SPr_{inj}(C^{f})$ est une cmf interne.

\begin{lem}\label{l2}
Le foncteur de localisation
$$Ho(SPr_{inj}(C^{f})) \longrightarrow Ho(L_{S}SPr_{inj}(C^{f}))$$
commute aux produits directs finis.
\end{lem}

\textit{Preuve du lemme:}
Par construction, le foncteur de localisation
est \'equivalent au foncteur
$$\mathbb{L}Re : Ho(SPr(C^{f})) \longrightarrow
Ho(M).$$
De plus, on a $\mathbb{L}Re(h_{x})\simeq
x$, et donc
$$\mathbb{L}Re(h_{x}\times h_{y})\simeq
\mathbb{L}Re(h_{x})\times^{h} \mathbb{L}Re(h_{y}),$$
pour tout $x,y\in C^{f}$.
Tout objet $F\in Ho(SPr(C^{f}))$ peut s'\'ecrire
comme une colimite homotopique d'objets de la forme
$h_{x}$ pour $x\in C^{f}$. On a donc, pour tout
$y\in C^{f}$
$$\mathbb{L}Re(F\times h_{y})\simeq
\mathbb{L}Re(Hocolim_{i\in I}h_{x_{i}}\times h_{y})\simeq
Hocolim_{i\in I}\left(\mathbb{L}Re(h_{x_{i}}\times h_{y})
\right) \simeq Hocolim_{i\in I}(x_{i}\times y))$$
$$\simeq
(Hocolim_{i\in I}x_{i})\times^{h} y\simeq
\mathbb{L}Re(Hocolim_{i\in I}h_{x_{i}})\times^{h} \mathbb{L}Re(F)(h_{y})
\simeq \mathbb{L}Re(F)\times^{h} \mathbb{L}Re(F)(h_{y}).$$
Finallement, pour deux objets
$F$ et $G$ de $Ho(SPr(C^{f}))$, on \'ecrit
$G\simeq Hocolim_{j}h_{y_{j}}$,
et on trouve
$$\mathbb{L}Re(F\times G)\simeq
Hocolim_{j}(\mathbb{L}Re(F\times h_{y_{j}}))\simeq
Hocolim_{j}(\mathbb{L}Re(F)\times^{h} y_{j})$$
$$\simeq \mathbb{L}Re(F)\times^{h} (Hocolim_{j}
\mathbb{L}Re(h_{y_{j}})) \simeq
\mathbb{L}Re(F)\times^{h} \mathbb{L}Re(G).$$
\hfill $\Box$ \\

Nous pouvons maintenant terminer la preuve de la
proposition. Nous savons que $SPr_{inj}(C^{f})$ est une cmf interne
et que le foncteur de localisation preserve les
produits directs (homotopiques).
On d\'eduit alors de fa\c{c}on purement formelle
que $L_{S}SPr_{inj}(C^{f})$ est encore une cmf interne. \hfill $\Box$ \\

On d\'eduit de Prop. \ref{p1} le corollaire important suivant.

\begin{cor}\label{cp1}
Soit $M$ une cmf combinatoire et faiblement
interne. La cat\'egorie $Ho(M)$ est cart\'esiennement close,
et ce de fa\c{c}on compatible avec le $Ho(SEns)$-enrichissement.
\end{cor}

Pour tout cmf combinatoire et faiblement interne $M$, nous
noterons $\mathbb{R}\underline{Hom}_{M}$ les $Hom$ internes
de la cat\'egorie homotopique $Ho(M)$. Ceci est bien entendu un
abus de notation, car $\mathbb{R}\underline{Hom}_{M}$ n'est pas
obtenu comme foncteur d\'eriv\'e d'un foncteur
$\underline{Hom}_{M}$ d\'efini sur $M$ elle m\^eme. 

Il est \`a noter que le corollaire \ref{cp1} peut \^etre
rafin\'e de la fa\c{c}on suivante. 
Soit $I$ une petite cat\'egorie, et soit 
$Y\in Ho(M^{I})$ et $X\in Ho(M)$. Alors, on peut d\'efinir
un objet $\mathbb{R}\underline{Hom}_{M}(X,Y)\in Ho(M^{I})$, 
fonctoriellement en $X$ et $Y$, de sorte \`a avoir 
des isomorphismes naturels dans $Ho(SEns)$
$$Map_{M^{I}}(Z\times^{h}X,Y) \longrightarrow
Map_{M^{I}}(Z,\mathbb{R}\underline{Hom}_{M}(X,Y)),$$
pour tout $Z\in Ho(M^{I})$. 
De m\^eme, on peut pour $X\in Ho(M^{I})$ et 
$Y\in Ho(M)$ d\'efinir un objet 
$\mathbb{R}\underline{Hom}_{M}(X,Y)\in Ho(M^{(I^{op})})$, 
fonctoriellement en $X$ et $Y$, et de sorte \`a avoir
des isomorphismes naturels dans $Ho(SEns)$
$$Map_{M^{(I^{op})}}(Z,\mathbb{R}\underline{Hom}_{M}(X,Y))\simeq
Hoend\left( (i,j)\mapsto Map_{M}(Z_{i}\times^{h}X_{j},Y) \right).$$ 
Par la suite, nous utiliserons
implicitement ces extensions functorielles des
$Hom$ internes de $Ho(M)$. 

\section{Intervalles}

Dans ce paragraphe on se fixe une
cmf combinatoire et faiblement interne $M$. \\

Pour un entier $n\geq 1$ et pour $0\geq i<n$, on dispose dans
$\Delta$ de morphismes
$$se_{i} : [1] \longrightarrow [n]$$
donn\'es en envoyant $0$ sur $i$ et
$1$ sur $i+1$. Ces morphismes v\'erifient
$se_{i}(1)=se_{i+1}(0)$, et induisent un isomorphisme
dans $\Delta$
$$\underbrace{[1]\coprod_{[0]}[1] \dots [1]\coprod_{[0]}[1]}_{n\; fois} 
\simeq [n].$$
Ainsi, pour tout objet
co-simplicial $C : \Delta \longrightarrow M$, et tout
entier $n\geq 1$ on dispose d'un morphism naturel dans $M$
$$C(1)\coprod_{C(0)}C(1) \dots
C(1)\coprod_{C(0)}C(1) \longrightarrow
C(n),$$
et donc d'un morphisme naturel dans $Ho(M)$
$$\underbrace{C(1)\coprod_{C(0)}^{\mathbb{L}}C(1) \dots
C(1)\coprod_{C(0)}^{\mathbb{L}}C(1)}_{n\; fois} \longrightarrow
C(n).$$

\begin{df}\label{d3}
Un objet co-simplicial
$$C : \Delta \longrightarrow M$$
est une \emph{co-cat\'egorie faible} si
pour tout $n>0$ le morphisme
$$C(1)\coprod_{C(0)}^{\mathbb{L}}C(1) \dots
C(1)\coprod_{C(0)}^{\mathbb{L}}C(1) \longrightarrow
C(n)$$
est un isomorphisme dans $Ho(M)$.
\end{df}

On se fixe pour la suite de ce paragraphe une co-cat\'egorie
faible
$$C : \Delta \longrightarrow M.$$

On consid\`ere les cmf $csM$ des objets co-simpliciaux
dans $M$, et $sM$ des objets simpliciaux
dans $M$ (pour la structure de Reedy par exemple, voir
\cite[\S 5.2]{ho}). 
Soit $X_{*} \in sM$ un objet simplicial
dans $M$. En utilisant la proposition \ref{p1} on voit que
l'on peut construire un foncteur bien d\'efini
$$\mathbb{R}\underline{Hom}_{M}(C,-) : Ho(M) \longrightarrow Ho(sM).$$

\begin{prop}\label{p3}
Avec les notations pr\'ec\'edentes, le foncteur $Ho(SEns)$-enrichi
$$\begin{array}{cccc}
|X_{*}|_{C} : & Ho(M)^{op} & \longrightarrow & Ho(SEns) \\
 & Y & \mapsto & Map_{sM}(X_{*},\mathbb{R}\underline{Hom}_{M}(C,Y))
\end{array}$$
est repr\'esentable par un objet $|X_{*}|_{C} \in Ho(M)$.
\end{prop}

\textit{Preuve:} \`A l'aide de la proposition \ref{p1} on voit que l'on peut
supposer que $M$ est une cmf simplicial et interne. On peut aussi supposer que
$C$ est cofibrant dans $csM$ et que 
$X_{*}$ est cofibrant dans $sM$ (pour les structures
de Reedy). Dans ce cas, 
on voit facilement que le
co-\'egaliseur des deux morphismes
naturels
$$ \coprod_{[m] \rightarrow [p]}C(m)\times X_{p} \rightrightarrows
\coprod_{[n]}C(n)\times X_{n} ,$$
repr\'esente le foncteur en question.
\hfill $\Box$ \\

\begin{df}\label{d4}
Pour un objet $X_{*}\in Ho(sM)$, l'objet
$|X_{*}|_{C}\in Ho(M)$ sera appel\'e
la \emph{$C$-r\'ealisation g\'eom\'etrique
de $X_{*}$}.
\end{df}

Bien entendu la construction $X_{*} \mapsto |X_{*}|_{C}$ est fonctorielle 
en $X_{*}$ et d\'efinit
un foncteur
$$|-|_{C} : Ho(sM) \longrightarrow Ho(M).$$

Soit $\widetilde{*}$ un mod\`ele cofibrant pour
l'objet final de $M$. On d\'efinit le foncteur
$$Ens \longrightarrow M$$
qui envoit un ensemble $E$ sur
$\coprod_{E}\widetilde{*}$. Ce foncteur s'\'etend naturellement
aux objets simpliciaux
$$\phi : SEns \longrightarrow sM.$$
Soit maintenant $A$ une petite cat\'egorie, et 
$N(A)\in SEns$ sont nerf. Par 
$\phi$ on en d\'eduit un objet
$\phi(N(A))\in Ho(sM)$. On pose alors
$$\overline{\Delta(1)}:=\phi(N(\overline{I})) \in Ho(sM),$$
o\`u $\overline{I}$ est la cat\'egorie avec 
deux objets et un unique isomorphisme entre eux.

\begin{df}\label{d5}
Une co-cat\'egorie faible $C$ dans $M$ est un
\emph{intervalle} si les deux propri\'et\'es suivantes
sont v\'erifi\'ees.
\begin{enumerate}
\item Le morphisme $C(0) \longrightarrow *$ est
une \'equivalence dans $M$.
\item Le morphisme $|\overline{\Delta(1)}|_{C} \longrightarrow *$
est une \'equivalence dans $M$.
\end{enumerate}
\end{df}

Supposons maintenant que
$C$ soit un intervalle au sens ci-dessus. 
Pour un diagrame 
$$\xymatrix{ & X \ar[d] \\
Y \ar[r] & Z}$$
dans $M$, on d\'efinit un objet
$X\times^{C}_{Z}Y\in Ho(M)$ par la formule
$$X\times^{C}_{Z}Y:=
(X\times^{h}Y)\times^{h}_{Z\times^{h}Z}\mathbb{R}\underline{Hom}_{M}(C(1),Z).$$
Dans cette formule, le morphisme
$$\mathbb{R}\underline{Hom}_{M}(C(1),Z) \longrightarrow 
Z\times^{h}Z$$
est induit par le morphisme naturel
$$*\coprod^{\mathbb{L}} *\simeq C(0)\coprod^{\mathbb{L}}C(0) \longrightarrow 
C(1).$$

\begin{df}\label{d6}
L'objet $X\times^{C}_{Z}Y\in Ho(M)$ est appel\'e le 
\emph{$C$-produit fibr\'e
homotopique de $X$ et $Y$ au-dessus de $Z$}.
\end{df}

On passera sous silence les propri\'et\'es de
fonctorialit\'e du $C$-produit fibr\'e homotopique, 
qui sont \'evidentes. \\

Pour terminer ce paragraphe, supposons 
maintenant que $p : X \longrightarrow Y$
soit un morphisme dans $M$. Nous allons d\'efinir
un objet simplicial $N^{C}(p)\in Ho(sM)$, appel\'e 
le $C$-nerf homotopique de $p$, et tel que
$$N^{C}(p)_{n}\simeq
\underbrace{X\times^{C}_{Y}X\times^{C}_{Y}\dots \times^{C}_{Y}X}_{n \; fois}.$$
Soit $Z_{*}\in Ho(sM)$, et consid\'erons
$|Z_{*}|_{C}$. On dipose d'un morphisme 
naturel $Z_{0}\times^{h}C(0) \longrightarrow |Z_{*}|_{C}$ dans $Ho(M)$, qui 
comme $C$ est un intervalle induit un morphisme naturel
$$Z_{0} \longrightarrow |Z_{*}|_{C}.$$

On peut donc consid\'erer
l'ensemble simplicial
$$Map_{M}(Z_{0},X)\times^{h}_{Map_{M}(Z_{0},Y)}
Map_{M}(|Z_{*}|_{C},Y).$$
Ceci d\'efinit un foncteur $Ho(SEns)$-enrichi
$$Z \mapsto Map_{M}(Z_{0},X)\times^{h}_{Map_{M}(Z_{0},Y)}
Map_{M}(|Z_{*}|_{C},Y)$$
de $Ho(sM)^{op}$ vers $Ho(SEns)$.

\begin{prop}\label{p4}
Le foncteur ci-dessus est repr\'esentable (au sens enrichi)
par un objet $N^{C}(p)\in Ho(sM)$. 
\end{prop}

\textit{Preuve:} On peut \`a l'aide de 
la proposition \ref{p1} supposer que $M$ est une cmf interne. 
De plus, on peut supposer que $C$ est un objet 
cofibrant dans $csM$ et que le morphisme 
$p$ est une fibration entre objets fibrants de $M$. 
On pose alors
$$N^{C}(p)_{n}:=
\underline{Hom}_{M}(C(0),X)^{n}\times_{\underline{Hom}_{M}(C(0),Y)^{n}}
\underline{Hom}_{M}(C(n),Y),$$
o\`u le morphisme 
$\underline{Hom}_{M}(C(n),Y)\longrightarrow\underline{Hom}_{M}(C(0),Y)^{n}$
est induit par les $n$-morphismes naturels 
$C(0) \longrightarrow C(n)$. Ceci d\'efinit 
l'objet $N^{C}(p)\in Ho(sM)$, et l'on remarque imm\'ediatement
qu'il poss\`ede la propri\'et\'e universelle requise. \hfill $\Box$ \\

\begin{df}\label{d7}
L'objet $N^{C}(p)\in Ho(sM)$ est appel\'e le 
\emph{$C$-nerf du morphisme $p$}.
\end{df}

Soit $Z_{*}\in Ho(sM)$ un objet simplicial. Alors il existe un 
morphisme naturel dans $Ho(sM)$
$$Z_{*} \longrightarrow N^{C}(p)$$
o\`u $p : Z_{0} \longrightarrow |Z_{*}|_{C}$ est le morphisme
naturel. Ceci se voit \`a l'aide de la propri\'et\'e universelle
du $C$-nerf.

\section{Th\'eorie de $(1,\infty)$-cat\'egories}

Supposons que $M$ soit une cmf faiblement interne, 
combinatoire, et $C\in Ho(csM)$ soit un intervalle. 
Nous aurons besoin des deux notions suivantes.

\begin{df}\label{d8}
\begin{enumerate}
\item
\emph{Une cat\'egorie faible dans $M$} est un 
objet simplicial $X_{*}\in sM$ tel que pour tout
$n\geq 1$ le morphisme
naturel
$$X_{n} \longrightarrow X_{1}\times^{h}_{X_{0}}X_{1}\dots
\times^{h}_{X_{0}}X_{1}$$
soit un isomorphisme dans $Ho(M)$ (en d'autres
termes c'est une cocat\'egorie faible dans
$M^{op}$).
\item Un objet $X\in M$ est \emph{$0$-local}
si le morphisme
$$Map_{M}(*,X) \longrightarrow Map_{M}(C(1),X)$$
est un isomorphisme dans $Ho(SEns)$.
\end{enumerate}
\end{df}

La d\'efinition fondamentale de ce travail est la suivante.

\begin{df}\label{dgiraud}
Une cmf $M$ munie d'un objet
$C\in Ho(csM)$ est \emph{une th\'eorie de $(1,\infty)$-cat\'egories}
si les conditions suivantes sont v\'erifi\'ees.
\begin{enumerate}
\item[A1] Pour tout objet fibrant et $0$-local $X$, la cmf 
$M/X$ est faiblement interne et 
combinatoire.
\item[A2] L'objet $C$ est un intervalle dans $M$.
\item[A3] Pour toute famille d'objets
$\{X_{i}\}$ dans $Ho(M)$, de somme $X=\coprod^{\mathbb{L}}X_{i}\in Ho(M)$,
et pour tout morphisme $Z\longrightarrow X$
dans $Ho(M)$, 
les morphismes naturels
$$X_{i} \longrightarrow X_{i}\times^{h}_{X}X_{i}\qquad
\coprod^{\mathbb{L}}Z\times^{h}_{X}X_{i} \longrightarrow Z$$
sont des isomorphismes dans $Ho(M)$ (pour tout indice $i$).
\item[A4] Soit $Z_{*}\in Ho(sM)$ une cat\'egorie faible, muni
d'un morphisme $|Z_{*}|_{C} \longrightarrow Y$ dans $Ho(M)$, et 
$p : X \longrightarrow Y$ un morphisme dans $M$, tel que 
$X$, $Z_{0}$ et $Z_{1}$ soient
$0$-locaux dans $M$.
Alors, le
morphisme naturel
$$|Z_{*}\times^{h}_{\mathbb{R}\underline{Hom}_{M}(C,Y)}X|_{C} \longrightarrow 
|Z_{*}|_{C}\times^{h}_{Y}X$$
est un isomorphisme dans $Ho(M)$. 
\item[A5] Pour toute cat\'egorie faible $X_{*} \in Ho(sM)$, telle que
les objets $X_{0}$ et $X_{1}$ soient $0$-locaux dans $M$,
le morphisme
naturel
$$X_{*} \longrightarrow N^{C}(X_{0} \rightarrow |X_{*}|_{C})$$
est un isomorphisme dans $Ho(sM)$. 
\item[A6] Un morphisme $f : X \longrightarrow Y$ dans
$M$ est une \'equivalence si et seulement si
les deux morphismes
$$Map_{M}(*,X) \longrightarrow 
Map_{M}(*,Y) \qquad Map_{M}(C(1),X) \longrightarrow 
Map_{M}(C(1),Y)$$
sont des isomorphismes dans $Ho(SEns)$.
\item[A7] Pour tout $[n]$ et $[m]$ dans 
$\Delta$ le morphisme
$$Hom_{\Delta}([n],[m]) \longrightarrow 
Map_{M}(C(n),C(m))$$
est un isomorphisme dans $Ho(SEns)$.
\end{enumerate}
\end{df}

Le lecteur remarquera l'analogie entre les axiomes
$A1-A5$ et les axiomes de Giraud caract\'erisant les
topoi ainsi que leur g\'en\'eralisations
aux topoi de mod\`eles de \cite{hagI}
(seul $A2$ est vraiment nouveau ici).
Les axiomes $A1$ et $A3$ nous disent essentiellement que
les colimites  homotopiques sont stables par 
produits fibr\'es homotopiques (du moins au-dessus
d'objets $0$-locaux), et que les sommes
sont homotopiquement disjointes. Les axiomes
$A4$ et $A5$ affirment quand \`a eux que les
relations d'\'equivalences \emph{cat\'egoriques}
sont effectives et universelles en un certains sens, la nouveaut\'e
ici \'etant bien de passer de la notion de relation d'\'equivalences
d\'efinie par des actions de groupoides \`a son 
analogue d\'efinie par des actions de cat\'egories. La condition 
$A6$ nous dit que $*$ et $C(1)$ sont g\'en\'erateurs, tout
comme on peut caract\'eriser la cat\'egorie des ensembles
comme le seul topos de Grothendieck 
tel que $*$ soit non vide et g\'en\'erateur (ici la condition
$A7$ sur l'intervalle $C$ implique automatiquement que $*$ et 
$C(1)$ sont non vides). La condition $A7$ 
est une fa\c{c}on d'exprimer que l'intervalle
est non d\'eg\'en\'er\'e, et en particulier non contractile.
Signalons enfin que du point de vue
intuitif $C(1)$ est la cat\'egorie avec deux objets
et un unique morphisme entre eux, et que de ce fait les
objets $0$-locaux correspondent aux 
$\infty$-cat\'egories qui sont des $\infty$-groupoides. \\

Il est important de remarquer qu'il existe au moins
une th\'eorie de $(1,\infty)$-cat\'egories. La plus
\'evidente est la th\'eorie homotopique des th\'eories
homotopiques de C. Rezk, d\'ecrite dans \cite{re}, que 
nous allons bri\'evement rappeler.

Notons $sSEns$ la cat\'egorie des objets simpliciaux
dans $SEns$, que l'on muni de sa structure de Reedy. Pour cette
structure les cofibrations et les \'equivalences faibles
sont d\'efinis niveaux par niveaux. Pour \'eviter des confusions
les objets de $sSEns$ seront toujours pens\'es comme des foncteurs
$$\begin{array}{cccc}
X_{*} : & \Delta^{op} & \longrightarrow & SEns \\
 & [n] & \mapsto & X_{n},
\end{array}$$
et non comme des ensembles bi-simpliciaux. 
On notera
$h(n)\in sSEns$ l'objet repr\'esent\'e par $[n]\in \Delta$, carat\'eris\'e
par l'isomorphisme d'adjonction 
$$Hom(h(n),X_{*})\simeq (X_{n})_{0}.$$ 
Dans $sSEns$, on dispose alors
des morphismes de Segal
$$\phi_{n} : h(1)\coprod_{h(0)}h(1) \dots h(1)\coprod_{h(0)}h(1) \longrightarrow
h(n),$$
pour tout $n\geq 2$. On dipose aussi du nerf 
$N(\overline{I})$, qui est un ensemble simplicial et vu
comme objet de $sSEns$ constant en la deuxi\`eme 
coordon\'ee ($\overline{I}$ d\'esigne toujours
la cat\'egorie avec deux objets et un unique
isomorphisme entre eux).
La cmf
$CSS$ est par d\'efinition la localis\'ee de Bousfield
\`a gauche de $sSEns$ le long des morphismes 
$N(\overline{I}) \rightarrow *$ et $\phi_{n}$ pour $n\geq 2$. 
Nous renvoyons \`a \cite{re} pour une description 
plus explicite et plus compl\`ete de la cmf
de $CSS$. Nous rappellerons
simplement la forme des objets fibrants. Un objet
$X_{*}\in sSEns$ est fibrant dans $CSS$ si c'est un 
espace de Segal complet qui est de plus Reedy fibrant. Cela
signifie que $X_{*}$ est fibrant pour la structure de 
Reedy sur $sSEns$, et de plus que les deux conditions suivantes
sont satisfaites.
\begin{itemize}
\item L'objet $X_{*}$ est un cat\'egorie faible dans $SEns$ au 
sens de la d\'efinition \ref{d8}. On d\'efinit
alors une cat\'egorie $[X_{*}]$ dont les objets 
sont les $0$-simplexes de $X_{0}$, et dont l'ensemble
de morphismes de $x$ vers $y$ est $\pi_{0}(Hom(x,y))$, 
o\`u $Hom(x,y)$ est la fibre homotopique
de $X_{1} \longrightarrow X_{0}\times X_{0}$ au point
$(x,y)$. La structure de cat\'egorie faible 
sur $X_{*}$ induit une structure naturelle de cat\'egorie
sur $[X_{*}]$.
\item On note $X_{hoequiv}$ le sous-ensemble simplicial
de $X_{1}$ form\'e de toutes les composantes
connexes $x\in \pi_{0}(X_{1})$ correspondant \`a 
des isomorphismes dans la cat\'egorie $[X_{*}]$. Alors 
le morphisme 
$$X_{0} \longrightarrow X_{1}$$
se factorise par 
$$X_{0} \longrightarrow X_{hoequiv}$$
et on demande que ce morphisme soit une \'equivalence
dans $SEns$. 
\end{itemize}

La premi\`ere des deux conditions pr\'ec\'edente
correspond \`a la localit\'e pour rapport 
aux morphismes $\phi_{n}$ pour $n>1$, et 
la seconde \`a la localit\'e par rapport
au morphisme $N(\overline{I}) \longrightarrow *$. De plus, 
lorsque $X$ est une cat\'egorie faible dans
$SEns$, il est important de rappeler que 
le morphisme naturel
$$Map_{sSEns}(N(\overline{I}),X_{*}) \longrightarrow X_{1}$$
induit une \'equivalence entre $Map_{sSEns}(N(\overline{I}),X_{*})$
est $X_{hoequiv}$. Nous utiliserons ceci \`a plusieurs
reprises par la suite.

La cmf $CSS$ dispose d'un intervalle naturel
$$C : \Delta \longrightarrow CSS$$ 
d\'efini
simplement par $C(n)=h(n)$. Le fait que $C$ soit un
intervalle provient pr\'ecisemment du proc\'ed\'e
de localisation passant de $sSEns$ \`a $CSS$ (noter
que $C$ n'est pas un intervalle dans $sSEns$). 
Noter que pour
cette th\'eorie, les objets $0$-locaux sont pr\'ecisement 
les objets simpliciaux $\Delta^{op} \longrightarrow SEns$
homotopiquement constants. 

On peut alors v\'erifier que la cmf $CSS$ munie
de son intervalle $C$ v\'erifient 
les axiomes de notre d\'efinition \ref{dgiraud}, et forment
donc une th\'eorie de $(1,\infty)$-cat\'egories. Nous ne le ferons
pas ici, et laissons le soin au lecteur de v\'erifier les 
axiomes $A1-A7$, qui ne pr\'esentent pas de grosses
difficult\'es et se d\'eduisent ais\'emment 
des r\'esultats de \cite{re}. 

\section{Th\'eor\`eme d'unicit\'e}

Notre th\'eor\`eme principal est le suivant. 

\begin{thm}\label{tgiraud}
Soit $(M,C)$ une th\'eorie de $(1,\infty)$-cat\'egories. 
Alors, $M$ est Quillen \'equivalente \`a $CSS$. 
\end{thm}

\textit{Preuve:} Par la proposition \ref{p1} on peut
supposer que $M$ est combinatoire, simpliciale
et que tous ces objets sont cofibrants. On supposera
aussi que $C\in csM$ est cofibrant pour la structure
de Reedy. Nous noterons $\underline{Hom}$ les
$Hom$ simpliciaux de $M$, $\underline{Hom}_{M}$
ses $Hom$ internes, et $X\otimes Z \in M$ le produit
tensoriel externe d'un ensemble simplicial $X$ par un objet 
$Z$ de $M$. On d\'efinit le foncteur
$S$ par la formule
$$\begin{array}{cccc}
S(X) : & \Delta^{op} & \longrightarrow & SEns \\
 & [n] & \mapsto & \underline{Hom}(C(n),X).
\end{array}$$
Ce foncteur poss\`ede clairement un adjoint 
\`a gauche $R$, qui envoit un 
objet $X_{*}\in CSS$ sur le co-\'egaliseur dans $M$ des
deux morphismes naturels
$$\coprod_{[m] \rightarrow [p]}X_{p}\otimes C(m) \rightrightarrows
\coprod_{[n]\in \Delta}X_{n}\otimes C(n).$$
Comme
$C$ est suppos\'e cofibrant dans $csM$, on voit 
facilement que l'adjonction $(R,S)$ est  une adjonction de Quillen
pour la structure de Reedy sur $sSEns$. 
Notons que $R(N(\overline{I}))$ est isomorphe
dans $M$ \`a $|\overline{\Delta(1)}|_{C}$. De m\^eme, 
l'image par $R$ du morphisme
$$\phi_{n} : h(1)\coprod_{h(0)}h(1) \dots h(1)\coprod_{h(0)}h(1) \longrightarrow
h(n),$$
est isomorphe au morphisme
$$C(1)\coprod_{C(0)}C(1) \dots C(1)\coprod_{C(0)}C(1) \longrightarrow
C(n).$$
Ainsi, comme
$C$ est un intervalle et par propri\'et\'e universelle
des localisation de Bousfield, l'adjonction de Quillen 
$(R,S)$ est aussi une adjonction de Quillen 
apr\`es localisation de $sSEns$ \`a $CSS$. Comme tout
objet de $CSS$ et de $M$ est cofibrant on a
clairement $\mathbb{L}R\simeq R$. 

Nous allons montrer que $(R,S)$ est une \'equivalence de Quillen.
Pour cela, notons que l'axiome $A6$ implique que le foncteur 
$$\mathbb{R}S : Ho(M) \longrightarrow Ho(CSS)$$
est conservatif. Il nous suffit donc de montrer que
pour tout $X_{*}\in CSS$, 
le morphisme d'adjonction 
$X_{*} \longrightarrow \mathbb{R}S\circ R(X_{*})$
est un isomorphisme dans $Ho(CSS)$. \\ 

Rappelons qu'un morphisme $f : X \longrightarrow Y$ dans
une cmf $N$ est un monomorphisme homotopique (ou simplement 
\emph{h-mono}) si 
le morphisme diagonal
$$X \longrightarrow X\times^{h}_{Y}X$$
est un isomorphisme dans $Ho(N)$. De mani\`ere \'equivalente, 
$f$ est un h-mono si et seulement si pour tout
objet $Z$ dans $N$ le morphisme
$$Map_{N}(Z,X) \longrightarrow Map_{N}(Z,Y)$$
est injectif sur les $\pi_{0}$ et bijectif sur
tous les $\pi_{i}$ (i.e. est \'equivalent \`a une
inclusion de composantes connexes).

\begin{lem}\label{l4}
Pour tout $X\in M$, le morphisme
$$\mathbb{R}\underline{Hom}_{M}(*,X) \longrightarrow
\mathbb{R}\underline{Hom}_{M}(C(1),X)$$
est un h-mono.
\end{lem}

\textit{Preuve du lemme \ref{l4}:} Par la propri\'et\'e
d'adjonction des $Hom$ internes il suffit de montrer que
pour tout $X\in M$ le morphisme
$$Map_{M}(*,X) \longrightarrow Map_{M}(C(1),X)$$
est un h-mono dans $SEns$. Mais par d\'efinition, ce
morphisme est \'equivalent au morphisme de 
d\'eg\'en\'erescence
$$\mathbb{R}S(X_{*})_{0} \longrightarrow \mathbb{R}S(X_{*})_{1}.$$
Comme $\mathbb{R}S(X_{*})$ est un objet fibrant de $CSS$, 
c'est un espace de Segal complet, et donc ce dernier 
morphisme est un h-mono par d\'efinition (cf \cite[\S 6]{re}). 
\hfill $\Box$ \\

\begin{lem}\label{l5}
Soit $X\in M$, et notons $\emptyset_{M} \in M$ l'objet
initial de $M$. S'il existe un morphisme
$X \longrightarrow \emptyset_{M}$ dans $Ho(M)$, alors
c'est un isomorphisme. De plus, les objets
$*$ et $\emptyset_{M}$ ne sont pas isomorphes dans $Ho(M)$.
\end{lem}

\textit{Preuve du lemme \ref{l5}:}
Commen\c{c}ons par montrer que 
$$Map_{M}(*,\emptyset_{M})=Map_{M}(C(1),\emptyset_{M})=\emptyset \in Ho(SSet),$$
o\`u $\emptyset$ d\'esigne l'objet initial de $SSet$.
Comme il existe un morphisme $* \longrightarrow C(1)$ dans
$Ho(M)$, on voit qu'il suffit de montrer que
$Map_{M}(*,\emptyset_{M})=\emptyset$. Pour cela, on
commence par utiliser l'axiome $A1$ qui implique que 
$X\times^{h}\emptyset_{M}\simeq \emptyset_{M}$ pour tout objet
$X\in M$. En prenant $X=C(1)$ on trouve que la projection
$$Map_{M}(*,\emptyset_{M})\times
Map_{M}(*,C(1)) \longrightarrow Map_{M}(*,\emptyset_{M})$$
est un isomorphisme dans $Ho(SEns)$. Enfin, l'axiome
$A7$ implique que 
$\pi_{0}(Map_{M}(*,C(1)))=*\coprod *$, ce qui implique
que $Map_{M}(*,\emptyset_{M})=\emptyset$. 

Supposons maintenant que $X \longrightarrow \emptyset_{M}$
soit un morphisme dans $Ho(M)$. On vient de voir que l'on doit avoir
$$Map_{M}(*,X)\simeq Map_{M}(C(1),X)\simeq 
Map_{M}(*,\emptyset_{M})\simeq Map_{M}(C(1),\emptyset_{M})\simeq
\emptyset.$$
L'axiome $A6$ implique donc que $X \longrightarrow \emptyset_{M}$
est un isomorphisme dans $Ho(M)$.

Finalement, le fait que $*$ et $\emptyset_{M}$ ne soient pas
isomorphes dans $Ho(M)$ se d\'eduit imm\'ediatement de l'axiome $A7$. 
 \hfill $\Box$ \\

\begin{lem}\label{l6}
Soit $X_{*} \in Ho(sM)$ une cat\'egorie faible
avec $X_{0}$ et $X_{1}$ des objets $0$-locaux dans $M$.
Alors, le morphisme naturel $X_{0} \longrightarrow |X_{*}|_{C}$ induit 
une surjection
$$[*,X_{0}] \longrightarrow [*,|X_{*}|_{C}].$$
\end{lem}

\textit{Preuve du lemme \ref{l6}:} Soit 
$* \longrightarrow |X_{*}|_{C}$ un morphisme dans
$Ho(M)$. On applique l'axiome $A4$ au diagrame
$$\xymatrix{
 & |X_{*}|_{C} \ar[d]^-{id} \\
\bullet \ar[r] & |X_{*}|_{C},}$$
qui nous donne un isomorphisme $|Y_{*}|_{C}\simeq *$ dans $Ho(M)$,
o\`u $Y_{*}$ est la cat\'egorie faible
d\'efinie par
$$Y_{*}:=X_{*}\times^{h}_{\mathbb{R}\underline{Hom}_{M}(C,|X_{*}|_{C})}*.$$ 
On a en particulier $Y_{0}\simeq X_{0}\times^{h}_{|X_{*}|_{C}}*$. 

Supposons que $[*,Y_{0}]=\emptyset$. Comme
$C$ est un intervalle, on a 
$[*,C(1)]\neq \emptyset$, et donc $[C(1),Y_{0}]=\emptyset$. 
Ainsi, par l'axiome $A6$ et le lemme \ref{l5} on en d\'eduit
que $Y_{0}\simeq \emptyset_{M}$. Le lemme \ref{l5} implique aussi
que $Y_{i}\simeq \emptyset_{M}$ pour tout $i$, et donc
que $|Y_{*}|_{C}\simeq *\simeq \emptyset$. Ceci est contradictoire
avec le lemme \ref{l5}, et ainsi on a 
$[*,Y_{0}]\neq \emptyset$. Comme la carr\'e suivant
$$\xymatrix{
Y_{0} \ar[r] \ar[d] & X_{0} \ar[d] \\
\bullet \ar[r] & |X_{*}|_{C}}$$
est homotopiquement cart\'esien, on en d\'eduit que le
morphisme $* \longrightarrow |X_{*}|_{C}$ se factorise
dans $Ho(M)$ par $X_{0} \longrightarrow |X_{*}|_{C}$, ce qu'il
fallait d\'emontrer.
\hfill $\Box$ \\

\begin{lem}\label{l7}
Soient $X_{*}$ et $Y_{*}$ deux cat\'egories faibles dans $M$ telles
que les deux conditions suivantes soient satisfaites.
\begin{enumerate}
\item L'objet $X_{0}$ est isomorphe dans $Ho(M)$ \`a
un objet de la forme $\coprod_{E}^{\mathbb{L}}*$, o\`u
$E$ est un ensemble. 
\item Les objets $Y_{0}$ et $Y_{1}$ sont $0$-locaux. 
\end{enumerate}
Alors, le morphisme naturel
$$[X_{*},Y_{*}]_{sM} \longrightarrow [|X_{*}|_{C},|Y_{*}|_{C}]_{M}$$
est surjectif. 
\end{lem}

\textit{Preuve du lemme \ref{l7}:} Soit 
$|X_{*}|_{C} \longrightarrow |Y_{*}|_{C}$ un morphisme
dans $Ho(M)$, et consid\'erons le diagramme suivant
$$\xymatrix{
X_{0}  \ar[d] & Y_{0} \ar[d] \\
|X_{*}|_{C} \ar[r] & |Y_{*}|_{C}.}$$
Le morphisme $X_{0} \longrightarrow |Y_{*}|_{C}$
est donn\'e par une application $E \longrightarrow
[*,|Y_{*}|_{C}]$. Ainsi, d'apr\`es le lemme \ref{l6}, 
il existe un morphisme $X_{0} \longrightarrow Y_{0}$ 
et un diagramme homotopiquement commutatif
$$\xymatrix{
X_{0}  \ar[r] \ar[d] & Y_{0} \ar[d] \\
|X_{*}|_{C} \ar[r] & |Y_{*}|_{C}.}$$
En passant aux $C$-nerfs, et en appliquant 
l'axiome $A5$ pour $Y_{*}$ on trouve un diagrame dans $Ho(sM)$
$$\xymatrix{
X_{*} \ar[r] \ar[rd] & N^{C}(X_{0}\rightarrow  |X_{*}|_{C}) \ar[r] \ar[d] & Y_{*} \ar[d] \\
 & |X_{*}|_{C} \ar[r] & |Y_{*}|_{C}.}$$
Ceci implique que le morphisme $X_{*} \longrightarrow Y_{*}$
induit par passage aux $C$-r\'ealisations g\'eom\'etriques
le morphisme original $|X_{*}|_{C} \rightarrow |Y_{*}|_{C}$. 
\hfill $\Box$ \\

\begin{lem}\label{l9}
Soit $L_{0}M$ la cmf localis\'ee de Bousfield
\`a gauche de $M$ le long du morphisme
$C(1) \longrightarrow *$. Alors, 
$L_{0}M$ est un topos de mod\`eles au sens
de \cite[Def. 3.8.1]{hagI}. 
\end{lem}

\textit{Preuve du lemme \ref{l9}:} On utilise le th\'eor\`eme
de Giraud pour les topos de mod\`eles
\cite[Thm. 4.9.2]{hagI}. 

\begin{sublem}\label{sl1}
Soit $X_{*}$ un groupoide
de Segal dans $M$ (au sens de \cite[Def. 4.9.1]{hagI}).
Alors, le morphisme
naturel 
$$|X_{*}|_{C} \longrightarrow |X_{*}|$$
est un isomorphisme dans $Ho(M)$.
\end{sublem}

\textit{Preuve du sous-lemme \ref{sl1}:}
Soit $X_{*}\in sL_{0}M$ un groupoide de Segal. 
On le repr\'esente par une cat\'egorie
faible $X_{*}$ dans $M$, tel que le morphisme
$$d_{0}\times d_{1} : X_{2} \longrightarrow X_{1}\times^{h}_{d_{0},X_{0},d_{0}}X_{1}$$
soit un isomorphisme dans $Ho(M)$.
Notons $|X_{*}|$ la colimite homotopique
de $X_{*}$ dans $M$, et consid\'erons le morphisme
naturel $|X_{*}|_{C} \longrightarrow |X_{*}|$. 

Soit $Z\in M$ un objet fibrant, et consid\'erons
le morphisme induit par $q$
$$q^{*} : Map_{M}(|X_{*}|,Z) \longrightarrow Map_{M}(|X_{*}|_{C},Z).$$
Il est isomorphe au morphisme naturel (ici 
$Z$ sera consid\'er\'e comme un objet simplicial constant)
$$Map_{sM}(X_{*},Z) \longrightarrow Map_{sM}(X_{*},Z^{C}),$$
que l'on sait \^etre un h-mono par le lemme \ref{l4}.
Pour voir que ce morphisme est un isomorphisme dans $Ho(SEns)$, 
on peut appliquer le lemme de Yoneda pour la cmf $sM$
(voir \cite[\S 4.2]{hagI}), et 
voir qu'il suffit de montrer que pour 
morphisme $u : X_{*} \longrightarrow Z^{C}$ dans $Ho(sM)$, et  
tout $Y\in M$, le morphisme induit par $u$
$$Map_{M}(Y,X_{*}) \longrightarrow Map_{M}(Y,Z^{C})$$
se factorise dans $Ho(sSEns)$ par 
le h-mono
$$Map_{M}(Y,Z) \longrightarrow Map_{M}(Y,Z^{C}).$$
Mais, comme $X_{*}$ est un groupoide de Segal, 
$Map_{M}(Y,X_{*})$ est un groupoide de Segal
dans $SEns$. De plus, $Map_{M}(Y,Z^{C})$
est un espace de Segal complet au sens de \cite{re}. 
Ainsi, $Map_{M}(Y,X_{*}) \longrightarrow Map_{M}(Y,Z^{C})$
se factorise par le sous-espace de Segal complet
des \'equivalences dans $Map_{M}(Y,Z^{C})$, qui n'est
autre que l'image de $Map_{M}(Y,Z) \longrightarrow Map_{M}(Y,Z^{C})$. 
Ceci finit de montrer que le morphisme $q : |X_{*}|_{C} \longrightarrow |X_{*}|$
est un isomorphisme dans $Ho(M)$.\hfill $\Box$ \\

\begin{sublem}\label{sl2}
Pour tout groupoide de Segal
$X_{*}$ dans $M$, tel que $X_{0}$ et $X_{1}$ soient
$0$-locaux, l'objet 
$|X_{*}|_{C}$ est $0$-local.
\end{sublem}

\textit{Preuve du sous-lemme \ref{sl2}:} 
Soit $f : C(1) \longrightarrow |X_{*}|_{C}$ un morphisme
dans $Ho(M)$. En utilisant le lemme \ref{l7}, on peut repr\'esenter
ce morphisme comme la $C$-r\'ealisation g\'eom\'etrique
d'un morphisme dans $Ho(sM)$
$$\Delta(1) \longrightarrow X_{*},$$
o\`u $\Delta(1)$ est l'objet simplicial de $M$ repr\'esent\'e
par $[1]$. Ce morphisme est aussi donn\'e par un morphisme
dans $Ho(sSEns)$
$$h(1) \longrightarrow Map_{M}(*,X_{*}).$$
Mais comme $X_{*}$ est un groupoide de Segal, $Map_{M}(*,X_{*})$ est un 
groupoide de Segal dans $SEns$, et on a donc
$Map_{M}(*,X_{*})_{hoequiv}=Map_{M}(*,X_{*})$. Ainsi, 
le th\'eor\`eme \cite[Thm 6.2]{re} implique que
le morphisme
pr\'ec\'edent se factorise en un diagramme commutatif
dans $Ho(sSEns)$
$$\xymatrix{
h(1) \ar[r]\ar[d] & Map_{M}(*,X_{*}) \\
N(\overline{I}). \ar[ru] & }$$
On en d\'eduit donc un diagramme commutatif dans $Ho(sM)$
$$\xymatrix{
\Delta(1) \ar[d] \ar[r] & X_{*} \\
\overline{\Delta(1)}, \ar[ru] & }$$
et donc un diagramme commutatif dans 
$Ho(M)$
$$\xymatrix{
C(1) \ar[d] \ar[r]^-{f} & |X_{*}|_{C} \\
|\overline{\Delta(1)}|_{C}. \ar[ru] & }$$
Comme $C$ est un intervalle dans $M$ on en d\'eduit que
le morphisme $f$ se factorise par $*$, et donc 
par le lemme \ref{l4} que le morphisme
$$Map_{M}(*,|X_{*}|_{C}) \longrightarrow Map_{M}(C(1),|X_{*}|_{C})$$
est un isomorphisme dans $Ho(SEns)$. Ceci implique
que $|X_{*}|_{C}$ est $0$-local.
\hfill $\Box$ \\

Revenons \`a la preuve du lemme \ref{l9}.
L'axiome $A1$ nous dit que 
la cmf $L_{0}M$ v\'erifie la condition
$1$ du th\'eor\`eme \cite[Thm. 4.9.2]{hagI}. 
La condition $2$ de \cite[Thm. 4.9.2]{hagI}
se d\'eduit de l'axiome A3. En effet, pour 
$\{X_{i}\}$ une famille d'objet dans
$Ho(L_{0}M)$, de somme $X=\coprod X_{i}$, et
$j$ un indice fix\'e,  le morphisme naturel 
$$\coprod_{i} X_{i}\times^{h}_{X}X_{j} \longrightarrow X_{j}$$
est un isomorphisme. Comme le morphisme
$$X_{j}\times^{h}_{X}X_{j} \longrightarrow X_{j}$$
est aussi un isomorphisme par $A3$, on voit que
pour tout $i\neq j$ on a
$$X_{i}\times^{h}_{X}X_{j}\simeq \emptyset_{L_{0}M}.$$
Ceci montre bien que $A3$ implique la condition $2$ de \cite[Thm. 4.9.2]{hagI}.

Pour terminer la preuve du lemma \ref{l9}, 
il ne nous reste donc qu'\`a d\'emontrer que
les relations d'\'equivalences de Segal
sont homotopiquement effectives dans $L_{0}M$ (au sens
de \cite[Def. 4.9.1]{hagI}). Soit $X_{*}$ un groupoide
de Segal dans $L_{0}M$ avec $X_{0}$ et $X_{1}$ $0$-locaux. 
D'apr\`es les sous-lemmes \ref{sl1} et \ref{sl2}, 
on sait que $|X_{*}|_{C}\simeq |X_{*}|$ est de plus que c'est un 
objet $0$-local. Le morphisme naturel 
$$X_{0}\times^{h}_{|X_{*}|}X_{0}\longrightarrow
X_{0}\times^{C}_{|X_{*}|_{C}}X_{0}$$
est donc un isomorphisme.
L'axiome $A5$ implique alors que le morphisme
$$X_{1} \longrightarrow X_{0}\times^{h}_{|X_{*}|}X_{0}$$
est un isomorphisme dans $Ho(L_{0}M)$, et donc que la relation
d'\'equivalence de Segal $X_{*}$ est homotopiquement effective.
\hfill $\Box$ \\

\begin{lem}\label{l10}
Soit $N$ un topos de mod\`eles. Suppons de plus que 
les trois propri\'et\'es suivantes soient satisfaites.
\begin{enumerate}
\item La cmf $N$ est simpliciale et tous ses objets
sont cofibrants.
\item Un morphisme $X \longrightarrow Y$ dans $N$
est une \'equivalence faible si et seulement si
le morphisme induit
$$Map_{N}(*,X) \longrightarrow Map_{N}(*,Y)$$
est un isomorphisme dans $Ho(SEns)$.
\item On a 
$$Map_{N}(*,\emptyset_{N})=\emptyset.$$
\end{enumerate}
Alors le foncteur de Quillen \`a droite
$$\underline{Hom}(*,-) : N \longrightarrow SEns$$
est une \'equivalence de Quillen.
\end{lem}

\textit{Preuve du lemme \ref{l10}:} La condition $3$ 
nous dit que le foncteur 
$$\mathbb{R}\underline{Hom}(*,-) : Ho(N) \longrightarrow
Ho(SEns)$$
est conservatif. Il nous reste donc \`a montrer que
pour tout $X\in SEns$, le morphisme d'adjonction
$$X\longrightarrow Map_{N}(*,X\otimes *)$$
est un isomorphisme dans $Ho(SEns)$. 

Tout d'abord, comme $N$ est un topos de mod\`eles, le
foncteur $X \mapsto X\otimes *$ commute aux produits fibr\'es 
homotopiques. De plus, la sous-cat\'egorie pleine $T$
de $Ho(N)$, form\'ee des objets $0$-tronqu\'es est un topos
de Grothendieck dans lequel l'objet final est g\'en\'erateur
et non vide. Le morphisme g\'eom\'etrique
$Ens \longrightarrow T$ est donc une \'equivalence de
cat\'egorie. En d'autres termes, pour tout ensemble
simplicial $0$-tronqu\'e $X$, le morphisme d'adjonction
$$X \longrightarrow Map_{N}(*,X\otimes *)$$
est un isomorphisme dans $Ho(SEns)$. 

Soit maintenant $X\in Ho(SEns)$, que l'on pr\'esente de
la forme $|Y_{*}|$ o\`u $Y_{*}$ est un groupoide de Segal
dans $SEns$ avec de plus $Y_{0}$ $0$-tronqu\'e.
On consid\`ere le morphisme 
$$Y_{0}\otimes * \longrightarrow X\otimes *\simeq |Y_{*}\otimes *|,$$ 
et le diagramme commutatif
$$\xymatrix{
Y_{0} \ar[r] \ar[d] & Map_{N}(*,Y_{0}\otimes *) \ar[d] \\
X \ar[r] &  Map_{N}(*,X\otimes *).}$$
Le morphisme horizontal du haut est un isomorphisme (car 
$T=Ens$). De plus, le morphisme $Y_{0}\otimes * \longrightarrow
X\otimes *$ est localement surjectif, et comme
$*$ ne poss\`ede pas de sous-objet non-trivial (par
les conditions $2$ et $3$) on en d\'eduit que le morphisme
vertical de droite est surjectif \`a homotopie pr\`es. 
Ceci implique clairement 
que le morphisme induit
$$\pi_{0}(X) \longrightarrow \pi_{0}(Map_{N}(*,X\otimes *))\simeq[*,X\otimes *]$$
est surjectif. Soit de plus $x$ et $y$ deux point de $Y_{0}$, 
tels que leurs images dans $[*,X\otimes *]$ soient \'egaux. Les deux
morphismes $x,y : * \longrightarrow Y_{0}$ sont donc \'egalis\'es
par $Y_{0} \longrightarrow X\otimes *$, et ils d\'efinissent donc
un morphisme dans $Ho(N)$
$$* \longrightarrow x\times^{h}_{X\otimes *}y\simeq (x\times^{h}_{X}y)\otimes *,$$
o\`u nous avons not\'e symboliquement $x\times^{h}_{X}y$
l'ensemble simplicial d\'efini par le diagramme
homotopiquement cart\'esien suivant
$$\xymatrix{
x\times_{X}^{h}y \ar[r]\ar[d] & \bullet \ar[d]^-{y} \\
\bullet \ar[r]_-{x} & X.}$$
Ainsi, $[*,(x\times^{h}_{X}y)\otimes *]\neq \emptyset$, et d'apr\`es
ce que l'on vient de voir ceci implique que 
$\pi_{0}(x\times^{h}_{X}y)$ est non vide, et donc
que $x=y$ dans $\pi_{0}(X)$. 

On a donc montr\'e que pour tout
$X\in Ho(SEns)$, le morphisme d'adjonction
$$\pi_{0}(X) \longrightarrow \pi_{0}(Map_{N}(*,X\otimes *))$$
est bijectif. Comme le foncteur
$X \mapsto X\otimes *$ commute de plus aux limites homotopiques finies,
ceci entraine aussi que le morphisme d'adjonction
$$X \longrightarrow Map_{N}(*,X\otimes *)$$
est un isomorphisme dans $Ho(SEns)$. 
\hfill $\Box$ \\

Un dernier lemme avant de revenir \`a la preuve
du th\'eor\`eme \ref{tgiraud}. 

\begin{lem}\label{l10+}
Soit $X$ un ensemble simplicial, alors
$X\otimes *$ est un objet $0$-local dans $M$.
\end{lem}

\textit{Preuve:} On commence par supposer que
$X$ est connexe. On peut alors \'ecrire, \`a \'equivalence
pr\`es, $X=BG$, o\`u $G$ est un groupe simplicial. 
Dans ce cas, il existe un isomorphisme naturel dans $Ho(M)$
$$X\otimes *\simeq |BG_{*}\otimes *|,$$ 
o\`u $BG_{*}\in sSEns$ est le classifiant de $G$
d\'efini comme d'habitude par $BG_{n}:=G^{n}$. 
L'objet simplicial $BG_{*}\otimes *\in sM$ est un groupoide 
de Segal dans $M$, et donc par le sous-lemme \ref{sl1}
on a $|BG_{*}\otimes *|_{C}\simeq |BG_{*}\otimes *|$. 
De plus, $BG_{0}\otimes *=*$ est un objet $0$-local. 

\begin{sublem}\label{sl3}
L'objet $BG_{1}\otimes *= G\otimes *$ est $0$-local.
\end{sublem}

\textit{Preuve du sous-lemme \ref{sl3}:} On consid\`ere
$\mathbb{R}S(BG_{*}\otimes *)$, qui est un groupoide
de Segal dans $CSS$, tel que 
$\mathbb{R}S(BG_{*}\otimes *)_{n}:=\mathbb{R}S(BG_{n}\otimes *)$.
En particulier, l'objet des objets de ce groupoide
de Segal est $\mathbb{R}S(BG_{0}\otimes *)\simeq *$.
Ainsi, la cat\'egorie homotopique 
$[\mathbb{R}S(BG_{*}\otimes *)_{1}]=[\mathbb{R}S(G\otimes *)]$
h\'erite d'une structure de groupoide de Segal
dans $Cat$ dont l'objet des objets est la
cat\'egorie ponctuelle $*$. Autrement dit, 
c'est une cat\'egorie monoidale dont la structure 
monoidale est inversible (\`a isomorphisme pr\`es). 
Il est bien connu que cela implique que 
la cat\'egorie $[\mathbb{R}S(G\otimes *)]$
est un groupoide, et donc que l'espace
de Segal complet $\mathbb{R}S(G\otimes *) \in CSS$
est homotopiquement constant. En d'autres termes 
l'objet $G\otimes *$ est $0$-local dans $M$. \hfill $\Box$ \\

Les sous-lemmes \ref{sl1}, \ref{sl2} et \ref{sl3} impliquent 
donc que l'objet $X\otimes *\simeq |BG_{*}\otimes *|_{C}$
est $0$-local dans $M$, lorsque $X$ est connexe. \\

Pour finir la preuve du lemme \ref{l10+} il nous
reste \`a traiter le cas g\'en\'eral
d'un ensemble simplicial $X\simeq \coprod X_{\alpha}$ avec
$X_{\alpha}$ connexe. On a bien entendu
$$X\otimes *\simeq \coprod X_{\alpha}\otimes *,$$
et par ce que l'on vient de voir il nous suffit de
d\'emontrer le sous-lemme suivant.

\begin{sublem}\label{sl4}
Soit $\{Y_{i}\}$ une famille d'objets $0$-locaux
dans $M$. Alors l'objet $Y:=\coprod Y_{i}$
est $0$-local.
\end{sublem}

\textit{Preuve du sous-lemme \ref{sl4}:}
Soit $C(1) \longrightarrow Y$ un morphisme 
dans $Ho(M)$. Il nous suffit de montrer qu'il 
existe un indice $j$ et un diagramme commutatif
dans $Ho(M)$
$$\xymatrix{
C(1) \ar[r] \ar[rd] & Y \\
 & Y_{j}. \ar[u]}$$
Par l'axiome $A3$ on a un isomorphisme dans $Ho(M)$
$$\coprod C(1)\times^{h}_{Y}Y_{i} \simeq C(1).$$
Cependant, en utilisant l'axiome $A7$ il est facile de
voir que $C(1)$ est ind\'ecomposable, et donc
qu'il existe un unique indice $j$ tel que
$$C(1)\times^{h}_{Y}Y_{j} \simeq C(1).$$
En effet, supposons qu'il existe deux indices
$i\neq j$, avec 
$$C(1)\times^{h}_{Y}Y_{i}\neq \emptyset_{M} \qquad
C(1)\times^{h}_{Y}Y_{j}\neq \emptyset_{M}.$$
Clairement, les deux morphismes distincts $* \longrightarrow C(1)$
(donn\'es par l'axiome $A7$)
se factorisent l'un par $C(1)\times^{h}_{Y}Y_{i}$ et
l'autre par $C(1)\times^{h}_{Y}Y_{j}$. 
On construirait alors un morphisme
$$C(1) \longrightarrow *\coprod *$$
qui envoit $C(1)\times^{h}_{Y}Y_{i}$ sur le premier point
et $C(1)\times^{h}_{Y}Y_{j}$ sur le second. En composant avec
l'automorphisme qui \'echange les deux points et par
le morphisme naturel
$$*\coprod * \longrightarrow C(1),$$ 
on trouverait un morphisme $C(1) \longrightarrow C(1)$
qui ne soit \'egal \`a aucun des trois endomorphismes
de $C(1)$ pr\'evus par l'axiome $A7$.
Il existe donc un unique indice $j$ tel que
$$C(1)\times^{h}_{Y}Y_{j} \simeq C(1),$$
ce qui implique l'existence de la factorisation souhait\'ee.
\hfill $\Box$ \\

On a ainsi termin\'e la preuve du lemme \ref{l10+}. \hfill $\Box$ \\

Nous revenons enfin \`a la preuve du th\'eor\`eme
\ref{tgiraud}. Soit $X_{*}\in CSS$ qui soit constant \'egal \`a $X_{0}\in SEns$.
Alors le sous lemme \ref{sl1} et le lemme \ref{l10+} impliquent
que 
$$R(X_{*})\simeq |X_{*}\otimes *|_{C}\simeq |X_{*}\otimes *|\simeq 
X_{0}\otimes *$$ 
est un objet $0$-local. Ainsi, les lemmes \ref{l9} et \ref{l10} impliquent que 
le
morphisme d'adjonction 
$$X_{*} \longrightarrow \mathbb{R}SR(X_{*}\otimes *)\simeq 
\mathbb{R}S(X_{*}\otimes*)$$
est un isomorphisme dans $Ho(CSS)$. La restriction du
foncteur 
$$R : Ho(CSS) \longrightarrow Ho(M)$$
\`a la sous-cat\'egorie des objets homotopiquement constants
est donc pleinement fid\`ele et son image consiste en tous
les objets $0$-locaux de $Ho(M)$. Par cette \'equivalence, on identifiera
donc $Ho(SEns)$ avec la sous-cat\'egorie pleine de $Ho(CSS)$ form\'ee
des objets homotopiquement constants, ainsi qu'avec la sous-cat\'egorie
pleine de $Ho(M)$ form\'ee des objets $0$-locaux. Pour \'eviter
des confusions nous continuerons \`a noter $X\otimes *$
l'image de $X\in Ho(SEns)$ dans $Ho(M)$ par cette identification.
Noter que l'on a
$$(X\times^{h}_{Z}Y)\otimes * \simeq (X\otimes *)
\times^{h}_{(Z\otimes *)}(Y\otimes *).$$

Soit $X_{*}\in Ho(CSS)$, et consid\'erons le morphisme naturel
$X_{0}\otimes * \longrightarrow R(X_{*})$. Nous avons 
$R(X_{*})\simeq |X_{*}\otimes *|_{C}$, et 
par le l'axiome $A5$, le morphisme naturel
$$X_{1}\otimes * \longrightarrow (X_{0}\otimes *)
\times^{C}_{R(X_{*})}(X_{0}\otimes *)$$
est un isomorphisme dans $Ho(M)$. 
On dispose d'un diagramme
commutatif dans $Ho(M)$
$$\xymatrix{
X_{0}\otimes * \ar[r] \ar[dr] & 
(X_{0}\otimes *)\times^{h}_{R(X_{*})}(X_{0}\otimes *) \ar[d]^-{a} \\
 & (X_{0}\otimes *)\times^{C}_{R(X_{*})}(X_{0}\otimes *)\simeq 
X_{1}\otimes *.}$$ 
D'apr\`es le lemme \ref{l4} le morphisme $a$ est un 
h-mono. De plus, comme $X_{*}$ est un espace de Segal
complet, le morphisme $X_{0}\otimes * \longrightarrow X_{1} \otimes *$
est aussi un h-mono. Ainsi, on en d\'eduit que 
$X_{0}\otimes * \longrightarrow (X_{0}\otimes *)\times^{h}_{R(X_{*})}
(X_{0}\otimes *)$
est un h-mono. 
De plus, il est facile de voir que l'image de
$(X_{0}\otimes *)\times^{h}_{R(X_{*})}(X_{0}\otimes *)$ dans
$X_{1}\otimes *$ est contenue dans le sous-objet 
$X_{hoequiv}\otimes *\subset X_{1}\otimes *$. 
Or, comme $X_{*}$ est un espace de Segal complet, on sait que
$X_{hoequiv}\otimes *$ est aussi l'image de $X_{0}\otimes * \longrightarrow 
X_{1}\otimes *$. 
Ceci implique que le morphisme
$$(X_{0}\otimes *) \longrightarrow (X_{0}\otimes *)\times^{h}_{R(X_{*})}
(X_{0}\otimes *)$$
est un isomorphisme dans $Ho(M)$. En d'autres termes, 
le morphisme $X_{0}\otimes * \longrightarrow R(X_{*})$ est un h-mono dans
$M$, et ainsi 
le morphisme
$$X_{0} \longrightarrow \mathbb{R}S(R(X_{*}))_{0}\simeq Map_{M}(*,R(X_{*}))$$
est un h-mono dans $SEns$. Comme le lemme \ref{l6} implique qu'il 
est surjectif \`a homotopie pr\`es, on voit que
c'est un isomorphisme dans $Ho(SEns)$. 

On consid\`ere maintenant le morphisme d'adjonction
$$X_{*} \longrightarrow \mathbb{R}SR(X_{*}).$$
C'est un morphisme entre espace de Segal complets, et l'on
vient de voir que c'est un \'equivalence en degr\'e $0$. 
De plus, le morphisme induit en degr\'e $1$
$$X_{1} \longrightarrow \mathbb{R}SR(X_{*})_{1}$$
est isomorphe dans $Ho(SEns)$ au morphisme naturel
$$X_{1}\simeq \mathbb{R}\underline{Hom}_{CSS}(h(1),X_{*})_{0} \longrightarrow
\mathbb{R}\underline{Hom}_{CSS}(h(1),\mathbb{R}SR(X_{*}))_{0}
\simeq Map_{M}(*,\mathbb{R}\underline{Hom}_{M}(C(1),R(X_{*}))).$$
Or, en utilisant ce que l'on vient de d\'emontrer 
$$Map_{M}(*,R(X_{*}))\simeq Map_{M}(*,X_{0}\otimes *),$$ 
nous trouvons
$$Map_{M}(*,\mathbb{R}\underline{Hom}_{M}(C(1),R(X_{*})))\simeq
Map_{M}(*,R(X_{*})\times^{C}_{R(X_{*})}R(X_{*}))$$
$$\simeq \left( Map_{M}(*,R(X_{*}))\times Map_{M}(*,R(X_{*})) \right) 
\times_{Map_{M}(*,R(X_{*})) \times Map_{M}(*,R(X_{*}))}
Map_{M}(*,\mathbb{R}\underline{Hom}_{M}(C(1),R(X_{*})))$$
$$\simeq 
\left( Map_{SEns}(*,X_{0})\times Map_{SEns}(*,X_{0}) \right) 
\times_{Map_{M}(*,R(X_{*})) \times Map_{M}(*,R(X_{*}))}
Map_{M}(*,\mathbb{R}\underline{Hom}_{M}(C(1),R(X_{*})))$$
$$\simeq Map_{M}(*,(X_{0}\otimes *)\times^{C}_{R(X_{*})}
(X_{0}\otimes *)).$$
Ainsi, le morphisme
$$X_{1} \longrightarrow \mathbb{R}SR(X_{*})_{1}$$
est isomorphe dans $Ho(SEns)$ au morphisme naturel
$$Map_{M}(*,X_{1}\otimes *) \longrightarrow Map_{M}(*,(X_{0}\otimes *)\times^{C}_{R(X_{*})}(X_{0}\otimes *)),$$
qui est un isomorphisme dans $Ho(SEns)$ par l'axiome
$A5$. Ceci finit de montrer que le morphisme d'adjonction 
$X_{*} \longrightarrow \mathbb{R}SR(X_{*})$ est un 
isomorphisme dans $Ho(CSS)$, et donc termine la preuve
du th\'eor\`eme 
\ref{tgiraud}.  \hfill $\Box$ \\

\section{Compl\'ements sur l'unict\'e}

Notre th\'eor\`eme \ref{tgiraud} affirme
que toute th\'eorie de $(1,\infty)$-cat\'egorie
est \'equivalente \`a $CSS$. Dans ce paragraphe nous
allons \'etudier l'unicit\'e de l'\'equivalence entre $CSS$ est une
th\'eorie de $(1,\infty)$-cat\'egorie, ou en d'autres termes
l'unicit\'e des auto-\'equivalences de $CSS$. \`A ce niveau, 
la th\'eorie des cmf n'est plus tr\`es bien adapt\'ee
car il est difficile de d\'efinir de fa\c{c}on pertinente
la notion d'auto-\'equivalences d'une cmf. Ce que nous
montrerons, c'est que la cat\'egorie simpliciale
$LCSS$, localis\'ee \`a la Dwyer et Kan de $CSS$ le long
des \'equivalences, poss\`edent un monoide
d'auto-\'equivalences \'equivalent \`a $\mathbb{Z}/2$.  

Pour cela, on rappelle l'existence d'une cmf des
$S$-cat\'egories, d\'emontr\'ee dans \cite{be}. 
La localis\'ee simpliciale $LCSS$ est donc un objet
de la cmf $S-Cat$ (le lecteur m\'eticuleux 
prendra garde de fixer
des univers). On consid\'ere
alors l'ensemble simplicial
$Map_{S-Cat}(LCSS,LCSS)$, ainsi que
$\mathbb{R}Aut(LCSS) \in Ho(SEns)$ le sous-ensemble simplicial
de $Map_{S-Cat}(LCSS,LCSS)$ form\'ee des 
\'equivalences. Le morphisme naturel 
$C : \Delta \longrightarrow LCSS$ permet de d\'efinir
$\mathbb{R}Aut(\Delta/LCSS)$ comme la fibre homotopique de
$$C^{*} : \mathbb{R}Aut(LCSS) \longrightarrow Map_{S-Cat}(\Delta,LCSS)$$
prise en le morphisme $C$. La proposition suivant nous
dit que l'objet $LCSS$ ne poss\`ede pas d'auto-\'equivalences
pr\'eservant le morphisme $\Delta \longrightarrow LCSS$.

\begin{prop}\label{pt2}
On a
$$\mathbb{R}Aut(\Delta/LCSS) \simeq *.$$
\end{prop} 

\textit{Esquisse de preuve:} Nous utiliserons l'analogue simplicial
du th\'eor\`eme principal de \cite[\S 7]{to}, dont nous
commencerons par rapeler l'\'enonc\'e. Supposons que 
$T$ et $T'$ soient deux $S$-cat\'egories.  
On forme les cmf $SPr(T)$ et $SPr(T')$ des
pr\'efaisceaux simpliciaux sur $T$ et $T'$ (munies
par exemples des structures projectives pour 
les quelles les fibrations et \'equivalences sont
d\'efinies termes \`a termes). On dispose alors
de deux nouvelles $S$-cat\'egories 
$$\widehat{T}:=Int(SPr(T)) \qquad
\widehat{T'}:=Int(SPr(T'))$$ 
form\'ees 
des objets fibrants et cofibrants 
dans $SPr(T)$ et $SPr(T')$. On consid\`ere
le sous ensemble simplicial 
$$Map_{S-Cat}^{c}(\widehat{T},\widehat{T'})
\subset Map_{S-Cat}(\widehat{T},\widehat{T})$$
form\'ee des morphismes qui commutent aux colimites
homotopiques\footnote{Nous passons sous silence
la notion g\'en\'erale de colimites homotopiques dans 
les $S$-cat\'egories, que l'on peut d\'efinir par exemple
\`a l'aide du plongement de Yoneda simplicial et de la
notion usuelle dans les cmf de pr\'efaisceaux simpliciaux.}. Alors, par exactement les m\^emes arguments
que ceux donn\'es pour les dg-cat\'egories dans \cite{to}
on montre que le plongement de Yoneda $l : T \longrightarrow \widehat{T}$
induit un isomorphisme dans $Ho(SEns)$.
$$l^{*} : Map_{S-Cat}^{c}(\widehat{T},\widehat{T'}) \longrightarrow
Map_{S-Cat}(T,\widehat{T'}).$$
Plus g\'en\'eralement, si $M$ est une cmf simpliciale et $LM\simeq Int(M)$ 
sa localis\'ee
simpliciale, alors le morphisme
$$l^{*} : Map_{S-Cat}^{c}(\widehat{T},LM) \longrightarrow
Map_{S-Cat}(T,LM).$$
est une \'equivalence.

Maintenant, soit $LCSS$ la localis\'ee simpliciale
de $CSS$ le long des \'equivalences. De m\^eme, soit 
$LSPr(\Delta)$ la localis\'e simpliciale de la cmf
des pr\'efaisceaux simpliciaux sur $\Delta$. 
L'adjonction de Quillen 
$$id : SPr(\Delta) \longrightarrow CSS \qquad 
SPr(\Delta) \longleftarrow CSS : id$$
induit des morphismes dans $Ho(S-Cat)$
$$l : LSPr(\Delta)\simeq \widehat{\Delta} \longrightarrow LCSS
\qquad \widehat{\Delta} \longleftarrow LCSS : i,$$
tels que $i$ soit pleinement fid\`ele (au sens
des $S$-cat\'egories) et avec $l\circ i=id$.  On d\'eduit de cela que
le morphisme
$$\xymatrix{
Map_{S-Cat}^{eq}(LCSS,LCSS) \ar[r]^-{l^{*}} & 
Map_{S-Cat}^{c}(\widehat{\Delta},LCSS)\simeq Map_{S-Cat}(\Delta,LCSS)}$$
admet une r\'etraction. Ainsi, le morphisme
$$\mathbb{R}Aut(\Delta/LCSS) \longrightarrow
Map_{\Delta/S-Cat}(\Delta,LCSS)\simeq *$$
admet aussi une r\'etraction. Ceci implique que
$\mathbb{R}Aut(\Delta/LCSS) \simeq *$, ce qu'il fallait 
d\'emontrer.
\hfill $\Box$ \\

\begin{prop}\label{pt22}
Tout auto-morphisme $f : LCSS \longrightarrow LCSS$
dans $Ho(S-Cat)$ pr\'eserve globalement 
l'image essentielle du foncteur $LC : \Delta \longrightarrow LCSS$. 
\end{prop}

\textit{Preuve:} On a clairement $f(*)\simeq *$. De plus, par la formule
$$C(n)\simeq C(1)\coprod_{*}^{\mathbb{L}} \dots C(1)\coprod_{*}^{\mathbb{L}}C(1)$$
on voit qu'il suffit de montrer que $f(C(1))$ est isomorphe dans
$Ho(LCSS)\simeq Ho(CSS)$ \`a $C(1)$. 

Pour cela, soit $Cat_{rig}$ la cat\'egorie des cat\'egories
rigides (i.e. dont les objets ne poss\`edent pas d'automorphismes
non triviaux). Le foncteur nerf induit un foncteur
pleinement fid\`ele $Cat_{rig} \longrightarrow Ho(CSS)$, et 
nous identifierons $Cat_{rig}$ \`a son image essentielle
par ce foncteur. 
Notons que 
le foncteur $C$ se factorise en 
$$C : \Delta \longrightarrow Cat_{rig} \longrightarrow Ho(CSS).$$
De plus, on peut caract\'eriser la sous-cat\'egorie
pleine $Cat_{rig}$ de $Ho(CSS)$ comme la sous-cat\'egorie
des objets $0$-tronqu\'es. Ainsi, l'auto-\'equivalence $f$ 
pr\'eserve forcemment la sous-cat\'egorie
$Cat_{rig}$. Pour d\'emontrer la proposition \ref{pt22}
il nous suffit donc de montrer que toute
auto-\'equivalence de la cat\'egorie 
$Cat_{rig}$ fixe, \`a \'equivalence pr\`es, la
cat\'egorie $I$ poss\`edant deux objets et un unique morphisme
entre eux. 

On remarque alors que $I$ est l'unique objet (\`a \'equivalence pr\`es)
poss\`edant les propri\'et\'es suivantes.
\begin{itemize}
\item On a $[*,I]_{CSS}\simeq *\coprod *$.
\item $I$ n'est pas isomorphe \`a $*\coprod *$.
\item Les seuls sous-objets de $I$ sont
$\emptyset$, $*$, $*\coprod *$ et $I$.
\end{itemize}
Comme ces propri\'et\'es sont invariantes par
toute auto-\'equivalence de $Cat_{rig}$, ceci termine la preuve
de la proposition. \hfill $\Box$ \\

\begin{thm}\label{t2}
Il existe une \'equivalence de monoides simpliciaux
$$\mathbb{R}Aut(LCSS)\simeq \mathbb{Z}/2.$$
\end{thm}

\textit{Preuve:} La proposition \ref{pt22} implique que toute auto-\'equivalence
de $LCSS$ induit une auto-\'equivalence de la cat\'egorie
$\Delta$. Ainsi, on dispose d'un morphisme de monoides simpliciaux
$$\mathbb{R}Aut(LCSS) \longrightarrow \mathbb{R}Aut(\Delta).$$
Il est facile de voir que l'on a
$\mathbb{R}Aut(\Delta)\simeq \mathbb{Z}/2$, la seule
auto-\'equivalence non triviale de la cat\'egorie $\Delta$ \'etant
celle qui fixe les objets $[n]$ en qui permute les deux
morphismes $[0] \rightrightarrows [1]$ (si on pense
\`a $\Delta$ comme la cat\'egorie des ensembles ordon\'es
finis, cette auto-\'equivalence est celle qui renverse
l'ordre). Notons $\sigma : \Delta \longrightarrow \Delta$
cet automorphisme. Il induit un automorphisme
de la cat\'egorie $CSS$, et donc de la $S$-cat\'egorie
$LCSS$. Ceci montre que le morphisme
$$\mathbb{R}Aut(LCSS) \longrightarrow \mathbb{R}Aut(\Delta)$$
est surjectif sur les $\pi_{0}$. De plus, la fibre homotopique
de ce morphisme est par d\'efinition $\mathbb{R}Aut(\Delta/LCSS)$,
qui est contractile par la proposition \ref{pt2}. Ceci finit
la preuve du th\'eor\`eme.  \hfill $\Box$ \\

\end{document}